\documentclass{article}
\usepackage[T1]{fontenc}
\usepackage{amsfonts}
\usepackage{amssymb}
\usepackage{amsmath}

\newtheorem{theorem}{Theorem}[section]
\newtheorem{lemma}[theorem]{Lemma}
\newtheorem{proposition}[theorem]{Proposition}
\newtheorem{corollary}[theorem]{Corollary}

\newenvironment{proof}[1][Proof]{\begin{trivlist}
\item[\hskip \labelsep {\bfseries #1}]}{\end{trivlist}}
\newenvironment{definition}[1][Definition]{\begin{trivlist}
\item[\hskip \labelsep {\bfseries #1}]}{\end{trivlist}}
\newenvironment{example}[1][Example]{\begin{trivlist}
\item[\hskip \labelsep {\bfseries #1}]}{\end{trivlist}}
\newenvironment{remark}[1][Remark]{\begin{trivlist}
\item[\hskip \labelsep {\bfseries #1}]}{\end{trivlist}}

\newcommand{\qed}{\nobreak \ifvmode \relax \else
      \ifdim\lastskip<1.5em \hskip-\lastskip
      \hskip1.5em plus0em minus0.5em \fi \nobreak
      \vrule height0.75em width0.5em depth0.25em\fi}

\begin{document}

\pagenumbering{arabic}

\title{Global Fluctuations in General $\beta$ Dyson Brownian Motion}
\author{Martin Bender}
\date{ }
\maketitle
\begin{abstract}
We consider a system of diffusing particles on the real line
 in a quadratic external potential and with repulsive 
electrostatic interaction. The empirical measure process 
is known to converge weakly to a deterministic measure-valued 
process as the number of particles tends to infinity. 
Provided the initial fluctuations are small, the rescaled 
linear statistics of the empirical measure process converge 
in distribution to a Gaussian limit for sufficiently smooth
 test functions. We derive explicit general formulae for the
 mean and covariance in this central limit theorem by analyzing
 a partial differential equation characterizing the limiting 
fluctuations.

   
\end{abstract}

\section{Introduction} We consider the following system of $n$ It\^o equations: 
\begin{equation}\label{SDE}d\lambda_t^i=\frac{2\sigma}{\sqrt{n\beta}}dB_t^i-\lambda_t^idt+\frac{2\sigma^2}{n}\sum_{j\neq i}\frac{dt}{\lambda_t^i-\lambda_t^j}, \textrm{ for $i=1,\ldots,n$}.
\end{equation} 
Here $\{B_t^i\}_{i=1}^n$ are independent, standard Brownian motions and $\sigma$ and $\beta>0$ are real parameters. These equations model the dynamics of $n$ diffusing particles on the real line with a logarithmic interaction potential, $u(x)=-\frac{1}{2}\log|x|$,  constrained by a quadratic external potential $v_n(x)=\frac{n x^2}{4 \sigma^2}$, at inverse temperature $\beta$. C\'epa and L\'epingle \cite{C&L} proved that the order of the particles is almost surely preserved for all times $t\geq0$. The stationary solution to (\ref{SDE}) has distribution 
\begin{multline}\frac{1}{\mathcal{Z}_n^{(\beta)}}\exp\left\{-\beta\left(\sum_{j=1}^nv_n(\lambda_j)+\sum_{i\neq j}u(\lambda_i-\lambda_j)\right)\right\}\prod_{i=1}^nd\lambda_i \\=\label{gibbs}\frac{1}{\mathcal{Z}_n^{(\beta)}}\prod_{1\leq i<j \leq n}|\lambda_i-\lambda_j|^{\beta}\exp\left\{-\frac{\beta n}{4\sigma^2}\sum_{j=1}^n\lambda_j^2\right\}\prod_{i=1}^nd\lambda_i,
\end{multline} 
where $\mathcal{Z}_n^{(\beta)}$ is a normalizing constant (the partition function) and $d\lambda$ denotes Lebesgue measure.

 For the specific parameter values $\beta=1,2 \textrm{ and } 4$ this model can also be interpreted in terms of matrix valued stochastic processes (Dyson Brownian motion). Let $\mathcal{M}_n(\beta)$ be the set of all $n\times n$ real ($\beta=1$), complex ($\beta=2$) and quaternion ($\beta=4$) matrices respectively and $\mathcal{S}_n(\beta)$ the set of self dual (with respect to conjugate transposition) elements in $\mathcal{M}_n(\beta)$. The Gaussian Orthogonal ($\beta=1$), Unitary ($\beta=2$) and  Symplectic ($\beta=4$) ensembles, $GX_{\beta}E_n(\sigma^2)$, with $X_{\beta}=O,U,S \textrm{ for } \beta=1,2,4$, are the probability distributions
\begin{equation}d\mu_n^{\beta}(M)=\frac{1}{\mathcal{Z}_n^{(\beta)}}\exp\left\{-\frac{\beta n}{4\sigma^2}\textrm{Tr} M^2\right\}dM \nonumber
\end{equation} 
on  $\mathcal{S}_n(\beta)$, where $dM=\prod_{i=1}^ndM_{ii}\prod_{1\leq i<j \leq n}dM_{ij}^{(1)}\cdots dM_{ij}^{(\beta)}$ is product Lebesgue measure on the essentially different members of $M=(M_{ij}^{(1)},\ldots,M_{ij}^{(\beta)})_{ij}$. 
 Let $M_{t}=(M^{ij}_{t})_{ij}$ be an  $\mathcal{S}_n(\beta)$-valued Ornstein-Uhlenbeck process, i.e. satisfying the SDE 
\begin{equation}\label{OU} dM_{t}=-M_{t}dt+\frac{\sigma}{\sqrt{\beta n}} d(B_{t}+B_{t}^{*}), \nonumber
\end{equation} 
where $B_{t}$ is an $n\times n$ matrix, the elements of which are independent standard real ($\beta=1$), complex ($\beta=2$) or quaternion ($\beta=4$) Brownian motions and $B_{t}^{*}$ is the conjugate transpose of $B_t$. Then the eigenvalues $\{\lambda_t^i\}_{i=1}^n$ of $M_{t}$ satisfy (\ref{SDE}) (see \cite{Ch}). For instance, if $M_0 \in GX_{\beta}E_n(\sigma_0^2)$  we will have $M_t \in GX_{\beta}E_n(e^{-2t}(\sigma_0^2-\sigma^2)+\sigma^2)$ for all $t\geq 0$ and if $M_0 \in\mathcal{S}_n(\beta)$ is fixed, equation (\ref{OU}) has solution $M_t=e^{-t}M_0+N_t,$ where $N_t\in GX_{\beta}E_n(\sigma^2(1-e^{-2t}))$. 

We define the empirical measure process \begin{equation} X_t^n=\frac{1}{n}\sum_{i=1}^n\delta_{\lambda_t^i}. \end{equation} To capture the asymptotic properties of the model on a global scale as $n\to \infty$, one is interested in studying the limiting behaviour of the linear statistics $\langle X_t^n, f \rangle =\frac{1}{n}\sum_{i=1}^nf(\lambda_t^i)$ where $f$ is a bounded continuous real test function.

\begin{example}\label{GUE+source} Define a \emph{deformed GUE} to be an ensemble of Hermitian matrices $M_a=M+D_a$ where $M$ is distributed according to the $GUE_{2n}(1)$ and $D_a=(d_{ij})_{i,j=1}^{2n}$ is a fixed $2n \times 2n$ diagonal matrix with 
\begin{equation}d_{ii}=\left\{\begin{array}{ll} a & \textrm{ for }1 \leq i \leq n \\ -a   &  \textrm{ for }n+1 \leq i \leq 2n\end{array}\right.\nonumber
.\end{equation} 
Then the eigenvalues of the rescaled matrix $M_a/\sqrt{a^2+\sigma^{-2}}$ correspond to the particles in our model with initial distribution $X_0=X_0^{2n}=\frac{1}{2}(\delta_{-1}+\delta_{1})$ at time $t=\log{\sqrt{1+(a\sigma)^{-2}}}$. The local behaviour of the eigenvalues in this model have been studied in \cite{B&H} and \cite{B&K} and it is known that the limiting eigenvalue density  of $M_a$ as $n\to \infty$ is supported on two disjoint intervals if $0<a<1$ and on one single interval if $a\geq1$. In other words,  $\textrm{supp}X_t$ grows from the two starting points $\{-1,1\}$ at time $t=0$ into two disjoint intervals that join at time  $t=\log{\sqrt{1+(a\sigma)^{-2}}}$.
\end{example}

In the stationary case, it is a classical result that  $\langle X^n, f \rangle$ converges in distribution to $\int_{\mathbb{R}} f d\mu,$ where 
\begin{equation}d\mu=\frac{1}{2\pi \sigma^2}\sqrt{4\sigma^2-x^2}\chi_{\{|x|<2\sigma\}}dx\nonumber
\end{equation}
 is the Wigner semi-circle law. 
More generally, for any initial asymptotic distribution of particles, $X_0$, the asymptotic particle distribution $X_t$ at each time $t\geq 0$ is uniquely determined by $X_0$ and converges weakly to $\mu$ as $t\to \infty$ (see Theorem \ref{Rogers&Shi} for a more precise statement).

A natural question is if there is a limiting distribution of the rescaled linear statistics. For ease of notation, we introduce the fluctuation process $Y_t^n=n(X_t^n-X_t)$, which takes signed Borel measures on $\mathbb{R}$ as values. We are interested in the limiting distribution of the random variables   
\begin{equation}\langle Y_t^n, f \rangle= \sum_{i=1}^nf(\lambda_t^i)-n\int_{\mathbb{R}} f(x) dX_t(x), 
\end{equation} 
where $f$ is a test function from an appropriate class, as $n\to \infty$. Note that there is no $\frac{1}{\sqrt{n}}$ normalization of the linear statistics here; this reflects the very regular spacing of the particles and is typical of related models.

 Israelsson \cite{Is} showed that $(\langle Y_{t_1}^n, f_1 \rangle,\ldots,\langle Y_{t_k}^n, f_k \rangle)$ has a Gaussian limit, provided the the test functions $f_j$ are $6$ times continuously differentiable and the initial distributions $X_0^n$ converge sufficiently fast to $X_0$. 
Although establishing existence and uniqueness, he does not characterize the limiting distribution-valued Gaussian process $Y_t$ very explicitly. In this work we derive explicit general formulae for the mean and covariance of the finite dimensional distributions of $Y_t$ by analyzing the partial differential equation arising in Israelsson's proof. 
 These formulae  generalize  many similar results obtained for various special cases of our model by completely different methods, some of these are briefly discussed below. In particular it is worth noting that our results hold for all values of the inverse temperature $\beta$ and in the non-equilibrium case with arbitrary initial particle distribution $X_0$.

Most of the previous related results pertain to specific matrix models and are restricted to the cases $\beta=1$ or $\beta=2$. The asymptotic global fluctuations for various ensembles of Hermitian and real symmetric matrices have been extensively studied, see e.g. \cite{B&S}, \cite{C}, \cite{G}, \cite{Jon}, \cite{P} and \cite{S&S}. In a recent paper \cite{B&Y}, Bai and Yao consider $N\times N$ matrices with zero mean, independent, not necessarily identically distributed entries such that the diagonal elements all have the same variance $\sigma^2/N$ and the off-diagonal elements have variance $1/N$ (real symmetric case) or uncorrelated real and imaginary parts each of variance $1/{2N}$ (Hermitian case). Such models are known as \emph{Wigner ensembles}. Under some fourth moment conditions, they provide a central limit theorem and give explicit mean and covariance formulae, which agree with those of Corollary (\ref{Equilibrium fluctuations}) with $\beta=2$ and $\Delta t=0$. Under the assumption of finite moments of all orders of all matrix elements, a more general class of ensembles of symmetric matrices is considered by Anderson and Zeitouni in \cite{A&Z}. Here the variances of all matrix elements, and the means of the diagonal entries, are allowed to depend on position. Spohn \cite{S} derives an expression for the covariance of the Gaussian fluctuations of our model in the hydrodynamic limit, but deals only with the case $\beta=2$ and (time dependent) equilibrium fluctuations. 

The few previous results available on the general $\beta$ case are restricted to an equilibrium situation. For the corresponding model on the circle, Spohn \cite{Sp} can handle the general $\beta$ case, again in the hydrodynamic limit at equilibrium. By expressing the equilibrium model in terms of ensembles of tridiagonal real matrices, Dumitriu and Edelman \cite{D&E} are able to find the general $\beta$ global fluctuations for polynomial test functions (this corresponds to the leading order term in Proposition \ref{polynomial}). Johansson \cite{J} considers a more general model corresponding to the equilibrium measure (\ref{gibbs}) but with the quadratic external potential $v_n$ replaced by a general polynomial of even degree and with positive leading coefficient. For the case of quadratic $v_n$ his mean and covariance formulae agree with the fixed $t$ equilibrium case of the model we discuss. 
In Johansson's model the variance is universal in the sense that it does not depend on the details of the potential, provided the support of the equilibrium measure is a single interval. In our model however, the variance at every finite $t$ depends on the initial conditions (see Proposition \ref{polynomial}). For instance, even though the eigenvalue density in the deformed GUE example will in finite time be supported on a single interval, the fluctuations remember the initial particle distribution for all $t\geq 0$. Thus the time evolution of the variance is determined by the geometry of the initial distribution; this structure is reminiscent of the role played by the boundary conditions in determining the fluctuations of the height function in discrete plane tiling models such as Kenyon's \cite{K}. There the  fluctuations converge to a Gaussian free field for a conformal structure determined by the boundary.

\paragraph{Acknowledgment}

I wish to thank Kurt Johansson for invaluable advice throughout the preparation of this work.

\section{Main results}
In order to formulate our results we need the following Theorem, referred to in the introduction: 
\begin{theorem}\label{Rogers&Shi}\emph{(Rogers and Shi, \cite{R&S}, C\'epa and L\'epingle, \cite{C&L})} Suppose that $X_0^n$ converges weakly in $\mathcal{M}$, the space of Borel probability measures on  $\mathbb{R}$ with the weak topology, to a point mass $X_0$ at an arbitrary element of $\mathcal{M}$.
Then there is a family $\{X_t\}_{t\geq 0} \subset \mathcal{M}$, depending only on $X_0$ and converging weakly as $t \to \infty$ to the Wigner semi-circle law, $\mu$, such that for each $t\geq 0$,  $X_t^n$ converges weakly to $X_t$ in $\mathcal{M}$ as $n \to \infty$. $X_t$ is uniquely characterized by the property that its Stieltjes transform, 
\begin{equation} M=M(t,z)=\int\frac{dX_t(x)}{x-z}, (t,z)\in [0,\infty)\times(\mathbb{C}\setminus \mathbb{R}), \nonumber 
\end{equation} 
solves the initial value problem 
\begin{equation}\label{R&S} \left\{\begin{array}{ll} M_t=(2\sigma^2M+z)M_z+M, & t>0\\M(0,z)=\int\frac{dX_0(x)}{x-z}.& {}\end{array}\right.
\end{equation}
\end{theorem}

We fix some terminology that will be used throughout the rest of this paper. Let $X_0$ be a given Borel probability measure on $\mathbb{R}$ and define $\Omega=\mathbb{C}\setminus\mathbb{R}$. Put $f(w)=\int\frac{dX_0(x)}{x-w}$, $w \in \Omega$; $f$ will be a holomorphic function. It follows from Theorem \ref{Rogers&Shi} that for every $t\geq 0$, $M(t,\cdot)=\int\frac{dX_t(x)}{x-\cdot}$ is a well-defined holomorphic function in $\Omega$, so we can define a family $\{h_t\}_{t\geq 0}$ of holomorphic maps in $\Omega$, 
\begin{equation}h_t(z)=ze^t+\sigma^2(e^t-e^{-t})M(t,z).
\end{equation}

\begin{proposition}\label{distribution} For every $t\geq 0$, $h_t(\Omega) \subseteq \Omega$ and the relation $g_t \circ h_{t}=\textrm{\emph{id}}$ holds, where 

\begin{equation} g_{t}(w)=e^{-t}w-\sigma^2(e^{t}-e^{-t})f(w).
\end{equation}  

Define $h_{t_1}^{t_2}=g_{t_2}\circ h_{t_1}$ for $t_1\geq t_2 \geq 0$. Then $h_{t_1}=h_{t_2}\circ h_{t_1}^{t_2}$.
\end{proposition}

\begin{proof} This is a step in the proof of Theorem \ref{chf}.
\qed
\end{proof}

Recall the definition of the Schwarzian derivative:
\begin{definition}\label{Schwartz} Let $v$ be a univalent function in some domain of the complex plane. The \emph{Schwarzian derivative} $Sv$ of $v$ is defined as 
\begin{equation} (Sv)(z)=\frac{v'''(z)}{v'(z)}-\frac{3}{2}\left(\frac{v''(z)}{v'(z)}\right)^2. \nonumber
\end{equation}  
We define the \emph{generalized Schwarzian derivative}, also denoted $Sv$, to be the following function of two complex variables: 
\begin{equation} (Sv)(z_1,z_2)=\left\{\begin{array}{ll}\frac{\partial^2}{\partial z_1 \partial z_2}\log\left(\frac{v(z_1)-v(z_2)}{z_1-z_2}\right)=\frac{v'(z_1)v'(z_2)}{(v(z_1)-v(z_2))^2}-\frac{1}{(z_1-z_2)^2} & \textrm{if $z_1\neq z_2$} \\ \lim_{z \to z_1}\frac{v'(z_1)v'(z)}{(v(z_1)-v(z))^2}-\frac{1}{(z_1-z)^2}=\frac{1}{6}(Sv)(z_1) &\textrm{if $z_1=z_2$}\end{array}\right.
\end{equation}
\end{definition}

\begin{definition} A probability measure on the space $\mathcal{S}'$ of tempered distributions is called an \emph{initial measure} if its law is a point mass at a $\nu\in \mathcal{S}'$ such that $\langle\nu,1\rangle=0$. 
\end{definition}

We can now state the main result, giving expressions for the mean and covariance of the finite dimensional distributions of the Stieltjes transform of the limiting fluctuation process $Y_t$.

\begin{theorem}\label{chf} \emph{(Mean and covariance formulae)} Let $Y_t^n=n(X_t^n-X_t)$ where $X_t^n$ is the empirical measure process and $X_t$ its weak limit. Suppose $Y_0^n$ converges weakly in $\mathcal{S}'$ to an initial measure $Y_0$ and that there is a constant $C$ such that for every $n$ and $z=a+bi$, $b\neq0$, the inequality 
\begin{equation}\mathbf{E}\left|\int \frac{d Y_0^n(x)}{x-z}\right|^2\leq\frac{C}{b^2} \nonumber
\end{equation} 
holds. Then $Y_t^n$ converges weakly to a Gaussian distribution-valued process $Y_t$ (see Theorem \ref{stefan} for the full statement), and for  $0\leq t_k\leq t_{k-1}\ldots \leq t_1$ and $z=(z_1,\ldots,z_k)\in (\mathbb{C}\setminus\mathbb{R})^k$ the Gaussian random vector 
\begin{equation}U=(U_1,\ldots,U_k) \textrm{, where } U_j=\langle Y_{t_j}, \frac{1}{\cdot-z_j}\rangle \nonumber 
\end{equation} 
has mean

\begin{equation} \mu_j=\frac{1}{2}\left(\frac{2}{\beta}-1\right)\frac{h_{t_j}''(z_j)}{h_{t_j}'(z_j)}+\langle Y_0,\frac{h'_{t_j}(z_j)}{\cdot-h_{t_j}(z_j)} \rangle 
\end{equation} 
and covariance matrix
\begin{multline}\Lambda_{lj}=\Lambda_{jl}=\frac{2}{\beta}\frac{\partial^2}{\partial z_j\partial z_l}\log \left(\frac{h_{t_j}(z_j)-h_{t_l}(z_l)}{h_{t_j}^{t_l}(z_j)-z_l}\right)\\
=\frac{2}{\beta}{h_{t_j}^{t_l}}'(z_j)\left(Sh_{t_l}\right)(h_{t_j}^{t_l}(z_j),z_l), \textrm{ if $l\geq j$}. 
\end{multline}

In particular,
\begin{equation}\label{varformula} \mathbf{Var}(U_j)=\frac{1}{3\beta}\left(Sh_{t_j}\right)(z_j).
\end{equation}
\end{theorem}
\begin{proof} See section \ref{bevis}.
\end{proof}

\begin{remark}An interesting consequence of equation (\ref{varformula}) is that a complete knowledge of the variance of the Stieltjes transform of $Y_t$ as a function of $z$ in the upper (or lower) half-plane at any fixed time $t\geq 0$, will uniquely determine the initial particle distribution $X_0$: By (\ref{varformula}) this function is the Schwarzian derivative of some analytic function, which is unique up to composition with an arbitrary M\"obius transformation (this is a well-known property of the Schwarzian derivative). It is easy to see that this determines $h_t$ uniquely, and this in turn determines the Stieltjes transform $f$ of the initial particle distribution $X_0$. 
\end{remark}
\begin{remark}
 Another model, concerning eigenvalues of non-Hermitian complex matrices, where a similar variance formula involving the Schwarzian derivative occurs is studied in \cite{W&Z}. 
\end{remark}

Let 
\begin{equation} f_{\mu}(z)=\frac{1}{2\pi \sigma^2}\int_{-2\sigma}^{2\sigma} \frac{\sqrt{4\sigma^2-x^2}}{x-z}dx=\frac{z}{2\sigma^2}\left(\sqrt{1-\left(\frac{2\sigma}{z}\right)^2}-1\right) \nonumber
 \end{equation} 
denote the Stieltjes transform of the Wigner semi-circle law, $\mu$. (Here $\sqrt{\cdot}$ means the branch of the square root for which $\Im(\sqrt{z})\geq 0$ iff $\Im z \geq 0$, defined for $z\in\mathbb{C}\setminus (-\infty,0)$.)

\begin{corollary}\label{Equilibrium fluctuations}\emph{(Equilibrium fluctuations) } Let $z_1, z_2\in \Omega$ and $\Delta t\geq 0$ be given. Put $t_1=t+\Delta t$ and $t_2=t$. Under the same hypotheses as in Theorem \ref{chf}, the asymptotic mean and covariance as $t\to \infty$ are given by
\begin{equation} \lim_{t\to \infty}m_j=\left(\frac{2}{\beta}-1\right)\frac{ \sigma^2 f_{\mu}(z_j)}{4\sigma^2-z_j^2},
\end{equation}
and
\begin{multline} \lim_{t\to \infty} \Lambda_{12}\\
=e^{-\Delta t}\frac{8\sigma^2 \left(\frac{1}{\sqrt{1-\left(\frac{2\sigma}{z_1}\right)^2}}+1\right)\left(\frac{1}{\sqrt{1-\left(\frac{2\sigma}{z_2}\right)^2}}+1\right)}{\beta \left(4\sigma^2 e^{-\Delta t}-z_1 z_2\left(\sqrt{1-\left(\frac{2\sigma}{z_1}\right)^2}+1\right)\left(\sqrt{1-\left(\frac{2\sigma}{z_2}\right)^2}+1\right)\right)^2}
\\=e^{-\Delta t}\frac{2\sigma^2  f_{\mu}'(z_1) f_{\mu}'(z_2)}{\beta (\sigma^2 f_{\mu}(z_1) f_{\mu}(z_2)e^{-\Delta t}-1)^2}. 
\end{multline}
\end{corollary}
\begin{proof} This is just a calculation using Theorem \ref{chf} and the fact that $X_t$ converges weakly to the semi-circle law (Theorem \ref{Rogers&Shi}).
\qed
\end{proof}

The previous results can be expressed in terms of integral formulae for the fluctuation process acting on analytic test functions. Suppose that for each $t\geq 0$ there is a compact set $C_t \subset \mathbb{R}$ such that $\mathrm{supp}X_t \subseteq C_t$ and, with probability $1$, $\mathrm{supp}Y_t \subseteq C_t$; we say that the fluctuation process $Y_t$ is compactly supported. Let $t_1\geq t_2\geq 0$ and for $i=1,2$ let  $\gamma_i$ be a closed simple curve in the complex plane, the interior of which contains $C_{t_i}$, and let $D_i$ be a simply connected domain containing $\gamma_i$.

\begin{theorem}\label{intrep} \emph{(Integral representation)} Suppose that $Y_t$ is compactly supported. Let $F_1$ and $F_2$ be analytic in the domains $D_1$ and $D_2$ defined above, respectively. Define the random variables $Z_1=\langle Y_{t_1},F_1\rangle$ and $Z_2=\langle Y_{t_2},F_2\rangle$. Then
\begin{equation}\label{intform}\mathbf{Cov}(Z_1,Z_2)=\frac{-1}{4\pi^2\beta}\oint_{\Gamma_1}\oint_{\Gamma_2}\left(F_1(g_{t_1}(w_1))-F_2(g_{t_2}(w_2))\right)^2\left(Sg_{t_2}\right)(w_1,w_2)dw_2dw_1, 
\end{equation}
where $\Gamma_i=h_{t_i}(\gamma_i)$. For $Z_1=Z_2$ this reduces further to 
\begin{equation}\mathbf{Var}(\langle Y_{t_1},F_1\rangle)=\frac{1}{4\pi^2\beta}\oint_{\Gamma_1}\label{intformvar}\oint_{\Gamma_1}\left(\frac{F_1(g_{t_1}(w_1))-F_1(g_{t_1}(w_2))}{w_1-w_2}\right)^2dw_2dw_1. 
\end{equation} 
\end{theorem}
\begin{remark} Note that the covariance depends only on the initial distribution $X_0$. \end{remark} 
\begin{proof} Representing $F_1$ and $F_2$ by the Cauchy integral formula as contour integrals along $\gamma_1$ and $\gamma_2$ we can use the linearity of $Y_t$ and Fubini's theorem to obtain
\begin{eqnarray}\label{cauchyint}
\mathbf{Cov}(Z_1,Z_2)&=&\mathbf{E}(Z_1Z_2)-\mathbf{E}Z_1\mathbf{E}Z_2 \nonumber\\
&=&\mathbf{E}\left[\langle Y_{t_1},\frac{1}{2\pi i}\oint_{\gamma_1}\frac{F_1(z)dz}{x-z}\rangle \langle Y_{t_2},\frac{1}{2\pi i}\oint_{\gamma_2}\frac{F_2(z)dz}{x-z}\rangle\right] \nonumber\\
&& -\mathbf{E}\left[\langle Y_{t_1},\frac{1}{2\pi i}\oint_{\gamma_1}\frac{F_1(z_1)dz}{x-z}\rangle \right] \mathbf{E}\left[ \langle Y_{t_2},\frac{1}{2\pi i}\oint_{\gamma_2}\frac{F_2(z)dz}{x-z}\rangle\right] \nonumber\\
&=&\oint_{\gamma_1}\oint_{\gamma_2}\frac{F_1(z_1)F_2(z_2)}{(2\pi i)^2}\left(\mathbf{E}\left[\langle Y_{t_1},\frac{1}{z_1-\cdot}\rangle \langle Y_{t_2},\frac{1}{z_2-\cdot}\rangle \right]\right.\nonumber \\
&&\left. -\mathbf{E}\left[\langle Y_{t_1},\frac{1}{z_1-\cdot}\rangle \right] \mathbf{E} \left[ \langle Y_{t_2},\frac{1}{z_2-\cdot} \rangle \right] \right) dz_2dz_1\nonumber \\
&=&\oint_{\gamma_1}\oint_{\gamma_2}\frac{F_1(z_1)F_2(z_2)}{(2\pi i)^2}\Lambda_{12}dz_2 dz_1,
\end{eqnarray}
where \begin{equation}\Lambda_{12}=\frac{2}{\beta}\frac{\partial^2}{\partial z_1\partial z_2}\log \left(\frac{h_{t_1}(z_1)-h_{t_2}(z_2)}{h_{t_1}^{t_2}(z_1)-z_2}\right) \nonumber
\end{equation} 
by  Theorem \ref{chf}. Since, for fixed $z_2$, $\Lambda_{12}$ is the derivative of an analytic function of $z_1$ in a domain containing $\gamma_1$, we note that
\begin{multline}\oint_{\gamma_1}\oint_{\gamma_2}(F_2(z_2))^2 \Lambda_{12}dz_2 dz_1\\
=\oint_{\gamma_2}(F_2(z_2))^2\oint_{\gamma_1}\frac{d}{dz_1}\left(\frac{-h_{t_2}'(z_2)}{h_{t_1}(z_1)-h_{t_2}(z_2)}+\frac{1}{h_{t_1}^{t_2}(z_1)-z_2}\right)dz_1dz_2=0.\nonumber
\end{multline}

Similarly,
\begin{equation}\oint_{\gamma_1}\oint_{\gamma_2}(F_1(z_1))^2 \Lambda_{12}dz_1 dz_2=0,\nonumber
\end{equation}
so we may substitute $-\frac{1}{2}\left(F_1(z_1)-F_2(z_2)\right)^2$ for the factor $F_1(z_1)F_2(z_2)$ in equation (\ref{cauchyint}), which  gives
 
\begin{multline}\mathbf{Cov}(Z_1,Z_2)=\frac{1}{4\pi^2\beta}\oint_{\Gamma_1}\oint_{\Gamma_2}\left(F_1(z_1)\right.\\
\left. -F_2(z_2)\right)^2\left( \frac{h_{t_1}'(z_1)h_{t_2}'(z_2)}{(h_{t_1}(z_1)-h_{t_2}(z_2))^2}-\frac{{h_{t_1}^{t_2}}'(z_1)}{(z_2-h_{t_1}^{t_2}(z_1))^2}\right) dz_2dz_1. \nonumber 
\end{multline} 
In the variance case, $Z_1=Z_2$, this reduces further since the second term of the integral becomes  
\begin{equation}\oint_{\Gamma_1}\oint_{\Gamma_2}\left(\frac{F_1(z_1)-F_1(z_2)}{z_2-z_1}\right)^2dz_2dz_1, \nonumber 
\end{equation} 
which vanishes by the analyticity of $F_1$. With the change of variables $w_1=h_{t_1}(z_1)$, $w_2=h_{t_2}(z_2)$ we arrive at the expression (\ref{intform}).
\qed
\end{proof}

As an application of Theorem \ref{intrep} we can consider polynomial test functions.
\begin{proposition} \emph{(Variance for polynomial test functions)}\label{polynomial} Suppose $Y_t$ is compactly supported. Then for $n=1,2, \ldots $ 
 \begin{equation} \label{polvar}\mathbf{Var}(\langle Y_t,x^{n}\rangle)=\left\{\begin{array}{ll} \frac{4\sigma^{2n}}{\beta}\sum_{s=1}^{n/2} s \binom{n}{n/2+s}^2+e^{-t}R_{n}(t,X_0) & \textrm{ if $n$ is even} \\e^{-t}R_{n}(t,X_0)  & \textrm{ if $n$ is odd,} \end{array}\right.
\end{equation} 
where $R_{k}$ is bounded in $t$ and depends only on moments of $X_0$ up to order $2k-1$.
\end{proposition}
\begin{proof}

Let $F_1(x)=F_2(x)=x^n$ and $t_1=t_2=t$.
After the change of variables $z_i=e^{t}/{w_i}$, formula (\ref{intformvar}) in this case reads 
\begin{equation}\label{variance}\mathbf{Var}(\langle Y_t,x^n\rangle)=\frac{1}{4\pi^2\beta}\oint_{|z_1|=r}\oint_{|z_2|=r}\left(\frac{g_t(e^{t}/z_1)^n-g_t(e^{t}/z_2)^n}{z_1-z_2}\right)^2dz_1dz_2, 
\end{equation}
where $r>0$ is such that $\textrm{supp} X_0\subset (-e^{t}/r,e^{t}/r)$. By definition of $g_t$, we can expand $g_t(e^{t}/z)$ in a Laurent series, 
\begin{equation}g_t(e^{t}/z)=\sum_{k=-1}^{\infty}a_k z^k, \nonumber
 \end{equation} where $a_{-1}=1$, $a_0=0$ and $a_k=\sigma^2(1-e^{-2t})e^{-(k-1)t}\int x^{k-1}dX_0(x)$ for $k\geq 1$, so equation (\ref{variance}) can be written 
\begin{multline}\label{variance2} \mathbf{Var}(\langle Y_t,x^n\rangle)\\
=\frac{1}{4\pi^2\beta} \oint_{|z_1|=r}\oint_{|z_2|=r}\left[\sum_{k_1,k_2,\dots,k_n\geq -1}a_{k_1}\cdots a_{k_n}\left(\frac{z_1^{k_1+\cdots +k_n}-z_2^{k_1+\cdots+k_n}}{z_1-z_2}\right)\right]^2dz_1dz_2. 
\end{multline}
For given integers $K$ and $J$ a simple combinatorial argument and the residue Theorem show that 
\begin{equation}\frac{1}{4\pi^2} \oint_{|z_1|=r}\oint_{|z_2|=r}\left(\frac{z_1^{K}-z_2^{K}}{z_1-z_2}\right)\left(\frac{z_1^{J}-z_2^{J}}{z_1-z_2}\right)dz_1dz_2= \left\{\begin{array}{ll} |K| & \textrm{if $K=-J$}\\
 0& \textrm{otherwise.}\end{array}\right.\nonumber 
\end{equation}
This means that we can write equation (\ref{variance2}) in the form
\begin{equation}\label{polysum}\mathbf{Var}(\langle Y_t,x^n\rangle)=\frac{2}{\beta}\sum_{s=1}^n s A_{-s,n}A_{s,n},
\end{equation}  
where 
\begin{equation} A_{s,n}=\sum_{\begin{array}{l}k_1+\cdots+k_n=s \\ k_i\geq -1 \end{array}} a_{k_1}\cdots a_{k_n}, \nonumber
\end{equation}
is a finite sum with the following structure: Since in the limit $t\to \infty$, $a_1\to \sigma^2$ and $a_k\to 0$ for $k>1$, all terms contributing to $A_{s,n}$ tend to $0$ exponentially in $t$ unless $k_i=\pm 1$ for $i=1,\ldots,n$. If $n+s$ is odd there are no such terms, and if $n+s$ is even there are $\binom{n}{\frac{n+s}{2}}$ choices of $(k_1,\ldots,k_n)$. Thus 
\begin{equation}\lim_{t\to \infty}A_{s,n}=\left\{\begin{array}{ll}0 & \textrm{ if $n+s$ is odd}\\ \binom{n}{\frac{n+s}{2}} \sigma^{n+s} & \textrm{ if $n+s$ is even.}\end{array}\right.\nonumber
\end{equation}
Inserting this into equation (\ref{polysum}) gives (\ref{polvar}).
\qed
\end{proof}
\begin{remark} The  $t \to \infty$ limit in this formula agrees with the variance formula of Dumitriu and Edelman \cite{D&E} and that of Johansson \cite{J} which  asserts that 
\begin{equation} \mathbf{Var}(\langle Y_t,h(x)\rangle)=\frac{1}{2\beta}\sum_{k=1}^{\infty}k\left(\frac{2}{\pi}\int_0^{\pi}h(2\sigma\cos(\Theta))\cos(k\Theta)d\Theta\right)^2 \nonumber
\end{equation} for an appropriate class of real test functions $h$.  Indeed, rewriting $h(x)=x^{2n}$ in terms of Chebyshev polynomials, the asymptotic variance is recovered from Johansson's result. 
\end{remark}

\section{Proof of Theorem 2.3}\label{bevis}
The proof of Theorem \ref{chf} relies on the characterization of $Y_t$ provided in \cite{Is} to prove existence and uniqueness of this process. For convenient reference we restate this result.

\begin{theorem}\label{stefan}\emph{(Israelsson,  \cite{Is})} Suppose that the sequence of measure-valued random variables $Y_0^n$ converges weakly in $\mathcal{S}'$ to an initial measure $Y_0$. Suppose further that there is a constant $C$ such that for every $n$ and $z=a+bi$, $b\neq0$, the inequality 
\begin{equation}\mathbf{E}\left|\int \frac{d Y_0^n(x)}{x-z}\right|^2\leq\frac{C}{b^2} \nonumber
\end{equation} holds.

Then $Y_t^n$ converges weakly to a Gaussian distribution-valued process $Y_t$ in the sense that for any $6$ times continuously differentiable, rapidly decreasing real test functions $f_j$, the random vector $\left(\int f_1(x) dY^n_{t_1}(x),\ldots,\int f_k(x) dY^n_{t_k}(x)\right)$ converges in distribution to $\left( \langle Y_{t_1},f_1\rangle,\ldots,\langle Y_{t_k},f_k\rangle \right)$.
Furthermore, the convergence extends to test functions of the form $\frac{1}{\cdot-z},$ $z\in\Omega=\mathbb{C}\setminus\mathbb{R}$, and $Y_t$ is uniquely characterized by its action on such functions by the following property: 
 Let  $0\leq t_{m+k}\leq t_{m+k-1}\ldots \leq t_{m+1}\leq t_1\leq T$ be given and for $s=(s_1,\ldots,s_m,\ldots,s_{m+k}) \in \mathbb{C}^{m+k}$, $z=(z_1,\ldots,z_m,\ldots,z_{m+k})\in \Omega^{m+k}$ and $t_{m+1}\leq t\leq t_1$  define the function 
\begin{multline}\phi(t,s_1,\ldots,s_m,z_1,\ldots,z_m)\\
=\mathbf{E}\left[\exp\left\{i\sum_{j=m+1}^{m+k} s_j\langle Y_{t_{j}},\frac{1}{\cdot-z_j} \rangle +i\sum_{j=1}^m s_j\langle{Y_t},\frac{1}{\cdot-z_j}\rangle \right\}\right]. \nonumber 
\end{multline} Then $\phi$ satisfies the PDE 

\begin{multline}\label{stefansPDE} \frac{\partial\phi}{\partial t}=\sum_{j=1}^m\left[s_j\left(1+2\sigma^2\frac{\partial M(t,z_j)}{\partial z_j}\right)\frac{\partial\phi}{\partial s_j}+\left(z_j+2\sigma^2 M(t,z_j)\right)\frac{\partial\phi}{\partial z_j}\right]\\
+\left\{2i\sigma^2\left(\frac{2}{\beta}-1\right)\sum_{j=1}^m\int\frac{s_j dX_t(x)}{(x-z_j)^3}-\frac{2\sigma^2}{\beta}\sum_{j=1}^m\sum_{l=1}^m\int\frac{s_j s_l dX_t(x)}{(x-z_j)^2(x-z_l)^2}\right\}\phi. \end{multline}

\end{theorem}
\begin{remark}
This is a slight reformulation of Israelsson's result: He allows for $Y_0$ to be random and works with real and imaginary parts of the complex functions $\frac{s}{\cdot-z}$ in order to ensure that the characteristic function $\phi$ be a priori well defined. However, once it is established that the distributions are Gaussian for such test functions if $Y_0$ is an initial measure, $\phi$ will be a well defined entire function of $s$ for test functions $\frac{1}{\cdot-z}$, $z\in\Omega$. The argument leading to equation (\ref{stefansPDE}) is then identical to that in Israelsson's proof, but this form is convenient for finding explicit solutions. 
\end{remark}  

\begin{remark}There is a numerical mistake in Israelsson's derivation of equation (\ref{stefansPDE}) which has been corrected here; all occurrences of the factor $\frac{\alpha}{2}$  in the equations on page $51$ and onward in \cite{Is} should be replaced by $\alpha$.
\end{remark}

Israelsson's method is similar to, although technically more involved than, that used by Rogers and Shi to prove Theorem \ref{Rogers&Shi}. The proof deals first with convergence for test functions in $\mathcal{S}$ and the convergence is then extended to the larger class of test functions. The existence of a limiting process $Y_t$ is shown by first establishing that the set of probability measures associated with the family $\{Y_t^n\}_{n=1}^{\infty}$ is tight in $C([0,T],\mathcal{S}')$. Then it is shown that any subsequential limit $\tilde{Y}_t$ must satisfy a certain martingale problem. By tightness, weak convergence of all finite dimensional distributions will entail the existence of a unique weak limit $Y_t$. This convergence is established by proving that any subsequential weak limit of the sequence $Y_t^n$ has the property that for every $f\in\mathcal{S}$ and $t\in[0,T]$ there is an approximating sequence $\{g_j\}_{j=1}^{\infty}$ of linear combinations of functions of the form $\frac{1}{\cdot-z}$ such that $\langle Y_t, g_j \rangle $ converges weakly to  $\langle Y_t,f \rangle $. This gives rise to the PDE (\ref{stefansPDE}) for the characteristic function of the Stieltjes transform of $\tilde{Y}_t$, which is shown to have a unique solution. Given an initial point mass distribution for $Y_0$, the resulting initial value problem is shown to have a unique solution which is the characteristic function of a Gaussian vector. Hence the finite dimensional distributions converge weakly by the approximation property and the theorem is proven. 

By solving equation (\ref{stefansPDE}) under appropriate initial conditions, we will be able to find expressions for the mean and covariance of the finite dimensional distributions of $Y_t$.
\begin{lemma}\label{IVP} For any fixed $t_0\geq 0$, let $\phi_{t_0}(s,z)$ be a given analytic function defined for $s=(s_1,\ldots,s_k) \in \mathbb{C}^k$ and $z=(z_1,\ldots,z_k)\in \Omega^k$ and let $U=\{(t,s,z)\lvert t>t_0, s\in \mathbb{C}^k, z\in \Omega^k\}$ and $\Gamma=\{(t_0,s,z)\lvert s\in \mathbb{C}^k, z\in \Omega^k\} \subseteq \partial U$. The initial value problem
\begin{multline} \frac{\partial\phi}{\partial t}=\sum_{j=1}^k\left[s_j\left(1+2\sigma^2\frac{\partial M(t,z_j)}{\partial z_j}\right)\frac{\partial\phi}{\partial s_j}+\left(z_j+2\sigma^2 M(t,z_j)\right)\frac{\partial\phi}{\partial z_j}\right]\\
 +\left\{2i\sigma^2\left(\frac{2}{\beta}-1\right)\sum_{j=1}^k\int\frac{s_j dX_t(x)}{(x-z_j)^3}-\frac{2\sigma^2}{\beta}\sum_{j=1}^k\sum_{l=1}^k\int\frac{s_j s_l dX_t(x)}{(x-z_j)^2(x-z_l)^2}\right\}\phi \textrm{ in $U$,} \nonumber
\end{multline}
\begin{equation}\label{PDE1}\phi(t_0,s,z)=\phi_{t_0}(s,z) \textrm{ on $\Gamma$,} \end{equation}

has the following unique solution:
\begin{equation} \label{solution}\phi(t,s,z)=\phi_{t_0}\left(s\cdot {h_t^{t_0}}'(z),h_t^{t_0}(z)\right)\exp\left\{i\sum_{j=1}^k s_j\mu_j-\frac{1}{2}\sum_{j=1}^k\sum_{l=1}^k s_j s_l\Lambda_{jl}\right\},\end{equation} where \begin{equation} \mu_j=\frac{1}{2}\left(\frac{2}{\beta}-1\right)\frac{{h_{t}^{t_0}}''(z_j)}{{h_{t}^{t_0}}'(z_j)},\end{equation}
 \begin{equation}\Lambda_{jl}=\Lambda_{lj}=\frac{2}{\beta}\left(Sh_{t}^{t_0}\right)(z_j,z_l), 
\end{equation}
and
 $s\cdot {h_t^{t_0}}'(z)$, $h_t^{t_0}(z)$ is shorthand notation for $(s_1 {h_{t}^{t_0}}'(z_1), \ldots s_k {h_{t}^{t_0}}'(z_k))$ and $ (h_{t}^{t_0}(z_1),\ldots,h_{t}^{t_0}(z_k))$ respectively.
\end{lemma}

\begin{proof}The equation is linear and can be solved with the method of characteristics. Fix $(t,s,z) \in U$ and let $\phi(\tau)=\phi(\overline{x}(\tau))$ be the solution along the characteristic $\overline{x}(\tau)=(t(\tau),s(\tau),z(\tau))$ through that point. By (\ref{PDE1}), the equations for $\overline{x}(\tau)$, if we choose $t(\tau)=\tau$, become
\begin{equation}\label{zeqn} \frac{dz_j(\tau)}{d\tau}=-z_j(\tau)-2\sigma^2M(\tau,z_j(\tau))
\end{equation} 
and
\begin{equation}\label{seqn} \frac{ds_j(\tau)}{d\tau}=-s_j(\tau)\left(1+2\sigma^2\frac{\partial M}{\partial z_j}(\tau,z_j(\tau))\right),
\end{equation}
while the solution $\phi(\tau)=\phi(\overline{x}(\tau))$ along the characteristic is given by the equation
\begin{multline}\label{phieqn} \frac{d\phi(\tau)}{d\tau}
=\left\{2i\sigma^2\left(\frac{2}{\beta}-1\right)\sum_{j=1}^k\int\frac{s_j(\tau) dX_t(x)}{(x-z_j(\tau))^3}\right.\\
\left. -\frac{2\sigma^2}{\beta}\sum_{j=1}^k\sum_{l=1}^k\int\frac{s_j(\tau) s_l(\tau) dX_t(x)}{(x-z_j(\tau))^2(x-z_l(\tau))^2}\right\}\phi(\tau).
\end{multline}
It will be convenient to solve these equations for all $\tau >0$ and impose the initial condition at the end. It may seem difficult to find solutions in closed form because of the dependence on the evolution of $X_t$, which is known only through the property of having Stieltjes transform satisfying (\ref{R&S}). As we will show however, all dependence on $X_t$ can be expressed in terms of $M$ and, more crucially, the evolution of $M$ along the characteristic is particularly simple. For the first point, it is easy to see by algebraic manipulations that
\begin{equation}\int\frac{dX_t(x)}{(x-z_j)^3}=\frac{1}{2}\frac{\partial^2}{\partial z_j^2}\left(\int\frac{ dX_t(x)}{(x-z_j)}\right)=\frac{1}{2}M_{zz}(t,z_j), \nonumber 
\end{equation}\begin{equation}\int\frac{dX_t(x)}{(x-z_j)^4}=\frac{1}{6}\frac{\partial^3}{\partial z_j^3}\left(\int\frac{ dX_t(x)}{(x-z_j)}\right)=\frac{1}{6}M_{zzz}(t,z_j), \nonumber 
\end{equation}

 and, with a little more effort,
\begin{equation}\int\frac{dX_t(x)}{(x-z_j)^2(x-z_l)^2}=\left(\frac{2(M(t,z_j)-M(t,z_l))}{(z_j-z_l)^3}-\frac{M_z(t,z_j)+M_z(t,z_l)}{(z_j-z_l)^2}\right),  \nonumber 
\end{equation} 
if $z_j\neq z_l.$ (Differentiating under the integral sign is clearly justified here since all integrands are bounded.)
Assuming without loss of generality that $z_j\neq z_l$ if $j\neq l$, equation (\ref{phieqn}) can thus be written 
\begin{multline}\label{phieqn2}\frac{1}{\phi(\tau)} \frac{d\phi(\tau)}{d\tau}\\
=2i\sigma^2\left(\frac{2}{\beta}-1\right)\sum_{j=1}^k\frac{s_j(\tau)}{2}M_{zz}(\tau,z_j(\tau))-\frac{2\sigma^2}{\beta}\left(\sum_{j=1}^k\frac{s_j(\tau)^2}{6}M_{zzz}(\tau,z_j(\tau)){} \right. \\
\left. +\sum_{j\neq l}s_j(\tau) s_l(\tau)\left(\frac{2(M(\tau,z_j(\tau))-M(\tau,z_l(\tau)))}{(z_j(\tau)-z_l(\tau))^3}-\frac{M_z(\tau,z_j(\tau))+M_z(\tau,z_l(\tau))}{(z_j(\tau)-z_l(\tau))^2}\right)\right)
\end{multline}
The equations (\ref{zeqn}) and (\ref{seqn}) can now be integrated with the aid of (\ref{R&S}) defining the evolution of $M(\tau,z_j(\tau))$. Fix $z_j=z$ and put $M(\tau)\equiv M(\tau,z(\tau))$, $M_z(\tau)\equiv \frac{\partial M(\tau,z(\tau))}{\partial z}$ and so on for all partial derivatives of $M(t,z)$. Differentiating, we have by the chain rule and equations (\ref{R&S}) and (\ref{zeqn}):
\begin{equation}\frac{dM(\tau)}{d\tau}=M_z(\tau)z'(\tau)+M_t(\tau)=M_z(\tau)(z'(\tau)+2\sigma^2M(\tau)+z(\tau))+M(\tau)=M(\tau).\nonumber 
\end{equation} 
or in integrated form simply 
\begin{equation}\label{M}M(\tau)=M(t,z)e^{\tau-t}.
\end{equation}
With (\ref{M}) substituted into (\ref{zeqn}), the latter equation can be integrated to yield 
\begin{equation} z(0)=ze^t+\sigma^2(e^{t}-e^{-t})M(t,z) \equiv h_t(z).\nonumber
\end{equation} 
Using this initial condition, equations (\ref{zeqn}) and (\ref{M}) give the explicit expression
\begin{equation}\label{ztau}z(\tau)=e^{-\tau}z(0)-\sigma^2(e^\tau-e^{-\tau})f(z(0))=g_{\tau}(z(0))=g_{\tau}(h_t(z))\end{equation} 
for the characteristic. In particular, taking $\tau=t$ gives $z=g_{t}(h_t(z))$, and since there is a unique characteristic through each point in $U$ it follows that for $t\geq t_1\geq t_0$, $h_{t}(z)=h_{t_1}(g_{t_1}(h_t(z)))$, which is the assertion of Proposition \ref{distribution}. Note that this provides a method of calculating $h_t$ (and $M(t,z)$) by finding an inverse of the explicitly defined function $g_t$.

Since the function $(t,z) \mapsto g_t(z)=e^{-t}z-\sigma^2(e^t-e^{-t})f(z)$ is $C^{\infty}$, it follows from implicit differentiaton of the relation $g_t(h_t(z))=z$ that the order of differentiation can be interchanged in the mixed partial derivatives of $h_t(z)$, in particular $\frac{\partial^{k+1}(h_t(z))}{\partial t \partial z^k}=\frac{\partial^{k+1}(h_t(z))}{\partial z^k\partial t}$ for $k=1,2,3$. Using this and differentating equation (\ref{R&S}) gives:
\begin{equation} M_{zt}=M_{tz}=(2\sigma^2M_z+1)M_z+(2\sigma^2M+z)M_{zz}+M_z, \nonumber 
\end{equation}
\begin{equation} M_{zzt}=M_{tzz}=(6\sigma^2M_z+3)M_{zz}+(2\sigma^2M+z)M_{zzz}, \nonumber 
\end{equation}
and
\begin{equation}M_{zzzt}=M_{tzzz}=6\sigma^2(M_{zz})^2+(8\sigma^2M_z+4)M_{zzz}+(2\sigma^2M+z)M_{zzzz}. \nonumber 
\end{equation} 
With these equations we can obtain ODE's for $M_z$, $M_{zz}$ and $M_{zzz}$ in a completely analogous fashion:
\begin{equation}\label{Mz}\frac{dM_z(\tau)}{d\tau}=M_{zz}(\tau)z'(\tau)+M_{zt}(\tau)=2(\sigma^2M_z(\tau)+1)M_z(\tau),
\end{equation}
\begin{equation}\label{Mzz}\frac{dM_{zz}(\tau)}{d\tau}=M_{zzz}(\tau)z'(\tau)+M_{zzt}(\tau)=(6\sigma^2M_z(\tau)+3)M_{zz}(\tau),
\end{equation}
and
\begin{equation}\label{Mzzz}\frac{dM_{zzz}(\tau)}{d\tau}=M_{zzzz}(\tau)z'(\tau)+M_{zzzt}(\tau)=6\sigma^2(M_{zz}(\tau))^2+4(2\sigma^2M_{z}(\tau)+1)M_{zzz}(\tau).
\end{equation}

Putting $w\equiv h_{t}(z)$, equation (\ref{M}) can be expressed

\begin{equation}\label{Mevolution}M(\tau)=f(w)e^{\tau},
\end{equation}
and equations (\ref{Mz}) through  (\ref{Mzzz}) can be integrated to produce 
\begin{equation}M_z(\tau)=f'(w)\frac{e^{\tau}}{g_{\tau}'(w)},
\end{equation}
\begin{equation}M_{zz}(\tau)=\frac{f''(w)}{(g_{\tau}'(w))^3},
\end{equation}
and
\begin{equation}M_{zzz}(\tau)=\left[f'''(w)+3(f''(w))^2\frac{\sigma^2(e^{\tau}-e^{-\tau})}{g_{\tau}'(w)}\right](g_{\tau}'(w))^{-4},
\end{equation} where $g_{\tau}(w)=e^{-\tau}w-\sigma^2(e^{\tau}-e^{-\tau})f(w)$. Inserting into equation (\ref{seqn}) and integrating we get 
\begin{equation}\label{stau}s(\tau)=s\frac{g_{\tau}'(w)}{g_{t}'(w)}. 
\end{equation} 

We can now finally express the right hand side of equation (\ref{phieqn2}) as an explicit function of $\tau$ by plugging in our expressions (\ref{ztau}) and (\ref{Mevolution}) through (\ref{stau}) derived for the evolution of $s_j$, $z_j$ and $z$-derivatives of $M$ along the characteristic. Integrating we see that 
\begin{equation}\label{logphi}\log\left(\frac{\phi(t,s,z)}{\phi(t_0,s(t_0),z(t_0))}\right)=I+II+III,
\end{equation} 
where
\begin{equation}I=2i\sigma^2\left(\frac{2}{\beta}-1\right)\int_{t_0}^t\sum_{j=1}^k\frac{s_j(\tau)}{2}M_{zz}(\tau,z_j(\tau))d\tau, 
\end{equation} 
\begin{equation} II=-\frac{2\sigma^2}{\beta}\int_{t_0}^t\left(\sum_{j=1}^k\frac{s_j(\tau)^2}{6}M_{zzz}(\tau,z_j(\tau))\right)d\tau, 
\end{equation}
 and 
\begin{multline}III=-\frac{2\sigma^2}{\beta}\int_{t_0}^t\sum_{j\neq l}s_j(\tau) s_l(\tau)\left(\frac{2(M(\tau,z_j(\tau))-M(\tau,z_l(\tau)))}{(z_j(\tau)-z_l(\tau))^3}\right. \\
\left. -\frac{M_z(\tau,z_j(\tau))+M_z(\tau,z_l(\tau))}{(z_j(\tau)-z_l(\tau))^2}\right)d\tau.
\end{multline}
To calculate these integrals we first note some immediate consequences of the definitions of the function $h_t^{t_0}=g_{t_0}\circ h_t$ and the generalized Schwarzian derivative:
\begin{equation}\label{secondderivative}
\frac{{h_{t}^{t_0}}''(z)}{{h_{t}^{t_0}}'(z)}=\frac{1}{g_t'(w)}\left(\frac{g_{t_0}''(w)}{g_{t_0}'(w)}-\frac{g_{t}''(w)}{g_{t}'(w)}\right),
\end{equation} 
and
\begin{equation}\label{Schwarzidentity} (Sh_t^{t_0})(z_1,z_2)=\frac{1}{g_t'(w_1)g_t'(w_2)}\left((Sg_{t_0})(w_1,w_2)-(Sg_{t})(w_1,w_2)\right), 
\end{equation}
where $w_i=h_t(z_i)$. Using the change of variables $x=b(\tau)=\sigma^2(e^{2\tau}-1)$ we can now calculate the integrals on the right hand side of equation (\ref{logphi}). First, we note that by (\ref{secondderivative}), 
\begin{multline}\int_{t_0}^t s(\tau) M_{zz}(\tau)d\tau
=s\int_{t_0}^t \frac{g_{\tau}'(w)}{g_{t}'(w)} \frac{f''(w)}{(g_{\tau}'(w))^3}d\tau\\
=\frac{sf''(w)}{2\sigma^2 g_t'(w)}\int_{b(t_0)}^{b(t)} \frac{dx}{(1-xf'(w))^2}
=\frac{sf''(w)}{2\sigma^2 g_t'(w) f'(w)}\left[\frac{1}{1-xf'(w)}-1\right]_{b(t_0)}^{b(t)}\\
=\frac{s}{2\sigma^2g_t'(w)}\left(\frac{e^{-t}b(t)f''(w)}{e^{-t}(1-b(t)f'(w))}-\frac{e^{-{t_0}}b(t_0)f''(w)}{e^{-{t_0}}(1-b(t_0)f'(w))}\right) =\frac{s}{2\sigma^2}\frac{{h_{t}^{t_0}}''(z)}{{h_{t}^{t_0}}'(z)}. \nonumber
\end{multline}
 This means that 
\begin{equation} I=2i\sigma^2\left(\frac{2}{\beta}-1\right)\int_{t_0}^t\sum_{j=1}^k\frac{s_j(\tau)}{2}M_{zz}(\tau,z_j(\tau))d\tau=\frac{i}{2}\left(\frac{2}{\beta}-1\right)\sum_{j=1}^k s_j\frac{{h_{t}^{t_0}}''(z_j)}{{h_{t}^{t_0}}'(z_j)}. 
\end{equation}
For the second integral,
\begin{eqnarray*}&&\int_{t_0}^t s(\tau)^2 M_{zzz}(\tau)d\tau \\&=&\frac{s^2}{(g_t'(w))^2}\int_{t_0}^t\frac{e^{2\tau}}{(1-b(\tau)f'(w))^2}\left[f'''(w)+\frac{3(f''(w))^2b(\tau)}{1-b(\tau)f'(w)}\right]d\tau{} \\
&=&\frac{s^2}{2\sigma^2(g_t'(w))^2}\int_{b(t_0)}^{b(t)}\frac{1}{(1-xf'(w))^2}\left(f'''(w)+\frac{3(f''(w))^2x}{(1-xf'(w))}\right)dx\\
&=&\frac{s^2}{2\sigma^2(g_t'(w))^2}\int_{b(t_0)}^{b(t)}\frac{1}{(1-xf'(w))^2}\left(\left(f'''(w)-\frac{3(f''(w))^2}{f'(w)}\right)\right.\\
&&\left. +\frac{3(f''(w))^2}{f'(w)}\frac{1}{(1-xf'(w))}\right) dx\\
&=&\frac{s^2}{2\sigma^2(g_t'(w))^2}\left[\left(\frac{f'''(w)}{f'(w)}-\frac{3(f''(w))^2}{(f'(w))^2}\right)\left(\frac{1}{(1-xf'(w))}-1\right)\right.\\
&&\left. +\frac{3}{2}\frac{(f''(w))^2}{(f'(w))^2}\left(\frac{1}{(1-xf'(w))^2}-1\right)\right]_{b(t_0)}^{b(t)} \\
&=&\frac{s^2}{2\sigma^2(g_t'(w))^2}\left(\frac{b(t)f'''(w)}{1-b(t)f'(w)}+\frac{3}{2}\left(\frac{b(t)f''(w)}{1-b(t)f'(w)}\right)^2 \right.\\
&&\left. -\frac{b(t_0)f'''(w)}{1-b(t_0)f'(w)}+\frac{3}{2}\left(\frac{b(t_0)f''(w)}{1-b(t_0)f'(w)}\right)^2\right) \\
&=&\frac{s^2}{2\sigma^2(g_t'(w))^2}\left(-(Sg_t)(w)+(Sg_{t_0})(w)\right)\\
&=&\frac{3s^2}{\sigma^2}(Sh_t^{t_0})(z,z), 
\end{eqnarray*} 
where we used the identity (\ref{Schwarzidentity}) in the last step.
Hence 
\begin{equation}II=-\frac{2\sigma^2}{\beta}\int_{t_0}^t\left(\sum_{j=1}^k\frac{s_j(\tau)^2}{6}M_{zzz}(\tau,z_j(\tau))\right)d\tau=-\frac{1}{\beta}\sum_{j=1}^k s_j^2(Sh_t^{t_0})(z_j,z_j).
\end{equation} 
 To calculate integral $III$, put $c=\frac{w_j-w_l}{f(w_j)-f(w_l)}$. Then for each $j\neq l$ we get a contribution to the sum in $III$ which takes the form
\begin{multline}\label{defint}
 -\frac{2\sigma^2}{\beta}\int_{t_0}^t s_j(\tau) s_l(\tau) \left( \frac{2(M(\tau,z_j(\tau))-M(\tau,z_l(\tau)))}{(z_j(\tau)-z_l(\tau))^3} \right.\\
\left.-\frac{M_z(\tau,z_j(\tau))+M_z(\tau,z_l(\tau))}{(z_j(\tau)-z_l(\tau))^2}\right) d \tau \\
= \frac{s_j s_l}{\beta g_t'(w_j)g_t'(w_l)\left(f(w_j)-f(w_l)\right)^2}\int_{b(t_0)}^{b(t)}\left( \frac{2cf'(w_j)f'(w_l)-(f'(w_j)+f'(w_l))}{(x-c)^2}\right.\\
\left. +\frac{2(cf'(w_j)-1)(cf'(w_l)-1)}{(x-c)^3} \right) dx.
\end{multline}
Now for any $s\geq0$ we can simplify 
\begin{multline}  \int_{0}^{b(s)}\left( \frac{2cf'(w_j)f'(w_l)-(f'(w_j)+f'(w_l))}{(x-c)^2} \right.\\
\left. +\frac{2(cf'(w_j)-1)(cf'(w_l)-1)}{(x-c)^3} \right) dx \nonumber
\end{multline}
\begin{eqnarray*}
&=& \left(2cf'(w_j)f'(w_l)-(f'(w_j)+f'(w_l))\right)\left(\frac{1}{b(s)-c}+\frac{1}{c}\right) \\
&&+\left(c^2f'(w_j)f'(w_l)-c(f'(w_j)+f'(w_l))+1\right)\left(\frac{1}{(b(s)-c)^2}-\frac{1}{c^2}\right)\\
&=&\frac{1}{c^2(b(s)-c)^2}\left(2cf'(w_j)f'(w_l)-(f'(w_j)+f'(w_l))\right)(cb(s)^2-c^2b(s))\\
&&+\left(c^2f'(w_j)f'(w_l)-c(f'(w_j)+f'(w_l))+1\right)(2cb(s)-b(s)^2)\\
&=&\frac{1}{c^2(b(s)-c)^2}\left(- (x-b(s))^2+c^2(b(s)f'(w_j)-1)(b(s)f'(w_l)-1)\right)\\
&=&(f(w_j)-f(w_l))^2\left(\frac{g_{s}'(w_j)g_{s}'(w_l)}{(g_{s}(w_j)-g_s (w_l))^2}-\frac{1}{(w_j-w_l)^2}\right)\\
&=&(f(w_j)-f(w_l))^2(Sg_s)(w_j,w_l),
\end{eqnarray*}
so by equations (\ref{defint}) and (\ref{Schwarzidentity}),
\begin{multline} III=\sum_{j\neq l}\frac{s_j s_l}{\beta g_t'(w_j)g_t'(w_l)}\left((Sg_t)(w_j,w_l)-(Sg_{t_0})(w_j,w_l)\right)\\
=-\frac{1}{\beta}\sum_{j\neq l}s_js_l(Sh_t^{t_0})(z_j,z_l).
\end{multline}
Inserting these integrals into equation  (\ref{logphi}) gives (\ref{solution}).\qed
\end{proof}

We are now ready to prove the main result, Theorem \ref{chf}. Let $s=(s_1,\ldots,s_k) \in \mathbb{C}^k$, $z=(z_1,\ldots,z_k)\in \Omega^k$ and $t=(t_1,\ldots,t_k)$, where  $0\leq t_k\leq t_{k-1}\ldots \leq t_1$, be given. We will prove that the characteristic function 
\begin{equation}\phi(t,s,z)=\mathbf{E}\left[\exp\left\{i\sum_{j=1}^{k} s_jU_j\right\}\right] \nonumber
\end{equation} 
of the random vector $U=(\langle Y_{t_1},\frac{1}{\cdot-z_1}\rangle, \ldots,\langle Y_{t_k},\frac{1}{\cdot-z_k}\rangle )$ is the characteristic function of a Gaussian vector. Since we have assumed that $Y_0$ is an initial measure it follows that $\phi_0(s,z)\equiv \phi(0,s,z)=\exp\left\{i\sum_{j=1}^{k} s_j \langle Y_{0},\frac{1}{\cdot-z_j}\rangle \right\}$. With the convention $t_{k+1}=0$, define the functions $\phi_{\tau}^{(j)}$, $j=1,\ldots, k$,  depending on the variables $s^{(j)}=(s_1^{(j)},\ldots,s_j^{(j)})$, $z^{(j)}=(z_1^{(j)},\ldots, z_j^{(j)})$ and the single time variable $\tau$, $t_{j+1}\leq \tau \leq t_j$ by the following expression: 
\begin{equation}\phi_{\tau}^{(j)}(s^{(j)},z^{(j)})=\mathbf{E}\left[\exp\left\{i\sum_{m=j+1}^k s_m\langle Y_{t_m},\frac{1}{\cdot-z_m}\rangle+i\sum_{m=1}^j s_m^{(j)}\langle Y_{\tau},\frac{1}{\cdot-z_m^{(j)}}\rangle\right\}\right]. \nonumber
\end{equation}  
Israelsson's Theorem (\ref{stefan}) states precisely that the $\phi_{\tau}^{(j)}$ satisfy equation (\ref{PDE1}) with initial conditions  $\phi_{t_{j+1}}^{(j)}(s^{(j)},z^{(j)})=\phi_{t_{j+1}}^{(j+1)}(s^{(j)},s_{j+1},z^{(j)},z_{j+1})$ for $j=1,\ldots,k-1$, and  $\phi_{0}^{(k)}(s^{(k)},z^{(k)})=\phi_{0}(s^{(k)},z^{(k)})$. Thus we may successively integrate $k$ times to obtain $\phi(t,s,z)=\phi_{t_1}^{(1)}(s_1,z_1)$ in terms of the initial conditions, using Lemma (\ref{IVP}) in each step. More explicitly, for $j=1,\ldots, k-1$ we have by  Lemma (\ref{IVP}) 
\begin{multline}\label{recursion} \phi_{t_j}^{(j)}(s^{(j)},z^{(j)})=\phi_{t_{j+1}}^{(j+1)}(s^{(j)}\cdot {h_{t_j}^{t_{j+1}}}'(z^{(j)}),s_{j+1}, {h_{t_j}^{t_{j+1}}}(z^{(j)}),z_{j+1})\\
\times \exp\left\{i\sum_{l=1}^j s_l^{(j)}\mu_l^{(j)}-\frac{1}{2}\sum_{l=1}^j\sum_{m=1}^j s_l^{(j)} s_m^{(j)}\Lambda_{lm}^{(j)}\right\}, 
\end{multline} 
where $\mu_l^{(j)}=\frac{1}{2}\left(\frac{2}{\beta}-1\right)\frac{{h_{t_j}^{t_{j+1}}}''(z_l^{(j)})}{{h_{t_j}^{t_{j+1}}}'(z_l^{(j)})} $ and $\Lambda_{lm}^{(j)}=\Lambda_{ml}^{(j)}=\frac{2}{\beta}\left(Sh_{t_j}^{t_{j+1}}\right)(z_l^{(j)},z_m^{(j)})$.

Applying formula (\ref{recursion}) $k-1$ times, starting with $\phi(t,s,z)=\phi_{t_1}^{(1)}(s_1,z_1)$, we obtain 
\begin{equation}\label{explicitphi}\phi(t,s,z)=\phi_{t_{k}}^{(k)}(s^{(k)},z^{(k)})\prod_{j=1}^{k-1}\exp\left\{i\sum_{l=1}^j s_l^{(j)}\mu_l^{(j)}-\frac{1}{2}\sum_{l=1}^j\sum_{m=1}^j s_l^{(j)} s_m^{(j)}\Lambda_{lm}^{(j)}\right\},
\end{equation}
where  $z_j^{(j)}=z_j$, $s_j^{(j)}=s_j$ and
\begin{equation} \label{recursivevariables}\left\{\begin{array}{ll} z_m^{(j)}=h_{t_{j-1}}^{t_{j}}\circ\ldots\circ h_{t_m}^{t_{m+1}}(z_m)=h_{t_m}^{t_{j}}(z_m) & \\ s_m^{(j)}=s_m\prod_{l=m}^{j-1} {h_{t_l}^{t_{l+1}}}'(z_m^{(l)})=s_m{h_{t_m}^{t_{j}}}'(z_m)  & \textrm{ for $m=1, \ldots, j-1$.} \end{array}\right.\nonumber
\end{equation} 
After a final application of Lemma (\ref{IVP}), equation (\ref{explicitphi}) becomes
\begin{multline} \label{explicitphi2}
\phi(t,s,z)=\phi_{0}(s\cdot h_t'(z),h_t(z))\exp\left\{\sum_{j=1}^k \left(i\sum_{l=1}^j s_l^{(j)}\mu_l^{(j)}-\frac{1}{2}\sum_{l=1}^j\sum_{m=1}^j s_l^{(j)} s_m^{(j)}\Lambda_{lm}^{(j)}\right)\right\}\\
=\phi_{0}(s\cdot h_t'(z),h_t(z))\exp\left\{i\sum_{l=1}^k\sum_{j=l}^k s_l^{(j)}\mu_l^{(j)}-\frac{1}{2}\sum_{l=1}^k\sum_{m=1}^k\sum_{j=max(l,m)}^k s_l^{(j)} s_m^{(j)}\Lambda_{lm}^{(j)}\right\}, 
\end{multline} 
where $s\cdot h_t'(z)=(s_1 h_{t_1}'(z_1), \ldots s_k h_{t_k}'(z_k))$ and  $h_t(z)=(h_{t_1}(z_1),\ldots,h_{t_k}(z_k))$.
 Now 
\begin{eqnarray*} s_l^{(j)}\mu_l^{(j)}&=& \frac{1}{2}\left(\frac{2}{\beta}-1\right)s_l{h_{t_l}^{t_j}}'(z_l)\frac{{h_{t_j}^{t_{j+1}}}''(h_{t_l}^{t_j}(z_l))}{{h_{t_j}^{t_{j+1}}}'(h_{t_l}^{t_j}(z_l))}\\
&=&\frac{1}{2}\left(\frac{2}{\beta}-1\right)s_l\frac{\frac{d}{dz_l}\left({h_{t_j}^{t_{j+1}}}'(h_{t_l}^{t_j}(z_l))\right)}{{h_{t_j}^{t_{j+1}}}'(h_{t_l}^{t_j}(z_l))}\\
&=&\frac{1}{2}\left(\frac{2}{\beta}-1\right)s_l\frac{d}{dz_l}\left(\log{h_{t_j}^{t_{j+1}}}'(h_{t_l}^{t_j}(z_l))\right), 
\end{eqnarray*} 
so since  
\begin{equation}\frac{d}{dz_l}\left(\log\prod_{j=l}^k{h_{t_j}^{t_{j+1}}}'(h_{t_l}^{t_j}(z_l))\right)
=\frac{d}{dz_l}\log ((h_{t_k}^{t_{k+1}}\circ\ldots\circ h_{t_l}^{t_{l+1}})'(z_l))
=\frac{d}{dz_l}\log(h_{t_l}'(z_l)),\nonumber
\end{equation}
we have found that 
\begin{equation}\sum_{j=l}^k s_l^{(j)}\mu_l^{(j)}=\frac{1}{2}\left(\frac{2}{\beta}-1\right)s_l\frac{d}{dz_l}\log(h_{t_l}'(z_l)).
\end{equation}
To evaluate $s_l s_m \Lambda_{lm}\equiv\sum_{j=max(l,m)}^k s_l^{(j)} s_m^{(j)}\Lambda_{lm}^{(j)}$ we distinguish between two cases. First, suppose that $h_{t_l}(z_l)\neq h_{t_m}(z_m)$. In this case, if we assume $l>m$, we can write 
\begin{eqnarray*}s_l s_m \Lambda_{lm}&=&\sum_{j=max(l,m)}^k s_l^{(j)} s_m^{(j)}\Lambda_{lm}^{(j)}\\
&=&s_l\sum_{j=l}^k{h_{t_l}^{t_j}}'(z_l)s_m{h_{t_m}^{t_j}}'(z_m) \frac{2}{\beta}\left(Sh_{t_j}^{t_{j+1}}\right)(h_{t_l}^{t_j}(z_l),h_{t_m}^{t_j}(z_m))\\
&=&s_l s_m\frac{2}{\beta}\sum_{j=l}^k{h_{t_l}^{t_j}}'(z_l){h_{t_m}^{t_j}}'(z_m)\\
&  &\times \frac{\partial^2}{\partial(h_{t_l}^{t_j}(z_l))\partial(h_{t_m}^{t_j}(z_m))}\log\left(\frac{h_{t_j}^{t_{j+1}}(h_{t_l}^{t_j}(z_l))-h_{t_j}^{t_{j+1}}(h_{t_m}^{t_j}(z_m))}{h_{t_l}^{t_j}(z_l)-h_{t_m}^{t_j}(z_m)}\right)\\
&=&s_l s_m\frac{2}{\beta}\sum_{j=l}^k \frac{\partial^2}{\partial z_l\partial z_m}\log\left(\frac{h_{t_j}^{t_{j+1}}(h_{t_l}^{t_j}(z_l))-h_{t_j}^{t_{j+1}}(h_{t_m}^{t_j}(z_m))}{h_{t_l}^{t_j}(z_l)-h_{t_m}^{t_j}(z_m)}\right)\\
& =&s_l s_m \frac{2}{\beta}\frac{\partial^2}{\partial z_l\partial z_m}\log\prod_{j=l}^k\left(\frac{h_{t_l}^{t_{j+1}}(z_l)-h_{t_m}^{t_{j+1}}(z_m)}{h_{t_l}^{t_j}(z_l)-h_{t_m}^{t_j}(z_m)}\right)\\
&=&s_l s_m \frac{2}{\beta}\frac{\partial^2}{\partial z_l\partial z_m}\log \left(\frac{h_{t_l}^{t_{k+1}}(z_l)-h_{t_m}^{t_{k+1}}(z_m)}{h_{t_l}^{t_l}(z_l)-h_{t_m}^{t_l}(z_m)}\right)\\
&=&s_l s_m \frac{2}{\beta}{h_{t_m}^{t_l}}'(z_m)\left(Sh_{t_l}\right)(z_l,h_{t_m}^{t_l}(z_m)),
\end{eqnarray*}
as claimed. Secondly, consider the case $h_{t_l}(z_l)= h_{t_m}(z_m)$. Using the identity $\left(S (f\circ g)\right)(z)=(g'(z))^2\left( Sf \right) (g(z))+\left(Sg\right)(z)$ for the Schwarzian derivative of a composition, we have 
\begin{eqnarray*}s_l s_m \Lambda_{lm}&=&\sum_{j=max(l,m)}^k s_l^{(j)} s_m^{(j)}\Lambda_{lm}^{(j)}\\
&=&\sum_{j=l}^k s_l{h_{t_l}^{t_j}}'(z_l)s_m{h_{t_m}^{t_j}}'(z_m) \frac{2}{\beta}\left(Sh_{t_j}^{t_{j+1}}\right)(h_{t_l}^{t_j}(z_l),h_{t_m}^{t_j}(z_m))\\
&=&s_l s_m \frac{2}{\beta}\sum_{j=l}^k{h_{t_l}^{t_j}}'(z_l){h_{t_m}^{t_j}}'(z_m) \frac{1}{6}\left(Sh_{t_j}^{t_{j+1}}\right)(h_{t_l}^{t_j}(z_l))\\
&=&s_l s_m\frac{2}{\beta}\sum_{j=l}^k{h_{t_l}^{t_j}}'(z_l){h_{t_m}^{t_j}}'(z_m)\frac{1}{6}\left(\frac{\left(S(h_{t_j}^{t_{j+1}}\circ h_{t_l}^{t_{j}}) \right)(z_l)-\left(Sh_{t_l}^{t_{j}}\right)(z_l)}{({h_{t_l}^{t_{j}}}'(z_l))^2}\right)\\
&=&s_l s_m\frac{2}{\beta}\sum_{j=l}^k{h_{t_m}^{t_l}}'(z_m)\frac{1}{6}\left(\left(Sh_{t_l}^{t_{j+1}} \right)(z_l)-\left(Sh_{t_l}^{t_{j}}\right)(z_l)\right)\\
&=&s_l s_m\frac{2}{\beta}{h_{t_m}^{t_l}}'(z_m)\left(Sh_{t_l}\right)(z_l,z_l).
\end{eqnarray*}
 Thus equation (\ref{explicitphi2}) can be expressed
\begin{multline}\phi(t,s,z)\\
=\exp\left\{i\sum_{j=1}^{k} s_j\left(h_{t_j}'(z_j)\langle Y_0,\frac{1}{\cdot-h_{t_j}(z_j)}\rangle +\frac{1}{2}\left(\frac{2}{\beta}-1\right)\frac{d}{dz_j}\log(h_{t_j}'(z_j))\right) \right. \\
\left.-\frac{1}{2}\sum_{l=1}^k\sum_{j=1}^ks_l s_j \frac{2}{\beta}{h_{t_j}^{t_l}}'(z_j)\left(Sh_{t_l}\right)(z_l,h_{t_j}^{t_l}(z_j))\right\},
\end{multline}
which shows that $U$ is Gaussian with mean and covariance as claimed.


\begin{thebibliography}{10}

\bibitem{A&Z} G. Anderson, O. Zeitouni, A CLT for a band matrix model, arXiv:math.PR/0412040 (2004)
\bibitem{B&K} P. Bleher, A. Kuijlaars, Large $n$ Limit of Gaussian Random Matrices with External Source, Part I, Communications in mathematical physics 252 (2004) 43-76
\bibitem{B&S} Z. Bai, J. Silverstein, CLT for linear spectral statistics of large-dimensional sample covariance matrices, Annals of Probability 32, no 1A, (2004) 553-605
\bibitem{B&Y} Z. Bai, J. Yao, On the convergence of the spectral empirical process of Wigner matrices, Bernoulli 11, no 6, (2005) 1059-1092
\bibitem{B&H} E. Br\'ezin, S. Hikami, Universal singularity at the closure of a gap in a random matrix theory, arXiv:cond-mat/9804023 (1998)
\bibitem{C} T. Cabanal-Duvillard, Fluctuations de la loi empirique de grande matrices al\'etoires, Annales de l'institut Henri Poincar\'e (B) Probability and statistics 37, no 3, (2001) 373-402 
\bibitem{C&L} E. C\'epa, D. L\'epingle, Diffusing particles with electrostatic repulsion, Probability Theory and Related Fields 107 (1997) 429-449 
\bibitem{Ch} T. Chan (1992), The Wigner semi-circle law and eigenvalues of matrix-valued diffusions, Probability Theory and Related Fields 93, no. 2, (1992) 249-272 
\bibitem{D&E} I. Dumitriu, A. Edelman, Global spectrum fluctuations for the $\beta$-Hermite and $\beta$-Laguerre ensembles via matrix models, arXiv:,ath-ph/0510043 (2006)
\bibitem{G} A. Guionnet, Large deviation upper bounds and central limit theorems for band matrices, Annales de l'institut Henri Poincar\'e (B) Probability and statistics 38, no 3, (2002) 253-384
\bibitem{Is} S. Israelsson, Asymptotic fluctuations of a particle system with singular interaction, Stochastic processes and their applications 93, no 1, (2001) 25-56
\bibitem{J} K. Johansson, On fluctuations of eigenvalues of random Hermitian matrices, Duke Mathematical Journal 91 (1998) 151-204
\bibitem{Jon} D. Jonsson, Some limit theorems for the eigenvalues of a sample covariance matrix, Journal of Multivariate analysis 12, no.1, (1982) 1-38
\bibitem{K} R. Kenyon, Height fluctuations in the honeycomb dimer model, arXiv.math-ph/0405052 (2004)

\bibitem{P} L. Pastur, Limiting laws of linear eigenvalue statistics for unitary invariant matrix models, preprint (2006)

\bibitem{R&S} L. Rogers, Z. Shi, Interacting Brownian particles and the Wigner law, Probability theory and related fields 95, no 4, (1993) 555-570

\bibitem{S&S} Y. Sinai, A. Soshnikov, Central limit theorem for traces of large random symmetric matrices with independent matrix elements, Boletim da Sociedade Brasileira de Matem\'atica 29, no.1, (1998) 1-24
\bibitem{S} H. Spohn, Interacting Brownian particles: a study of Dyson's model, in: Hydrodynamic behaviour and interacting particle systems, IMA Vol. Math. Appl.,9, 1986, pp.151-179
\bibitem{Sp} H. Spohn, Dyson's model of interacting Brownian motions at arbitrary coupling strength, Markov processes and related fields 4, no. 4, (1998) 649-661
\bibitem{W&Z} P. Wiegmann, A. Zabrodin, Large scale correlations in normal non-Hermitian matrix ensembles, Journal of physics A: Mathematical and general 36 (2003) 3411-3424
\end{thebibliography}
\end{document}